\documentclass[12pt]{amsart}
\usepackage{graphicx}
\usepackage{graphics}
\usepackage{amsmath}
\usepackage{amscd}
\usepackage{latexsym}
\usepackage{subfigure}
\usepackage{caption}

\begin{document}

\textwidth 5.9in
\textheight 7.9in

\evensidemargin .75in
\oddsidemargin .75in

\newtheorem{Thm}{Theorem}
\newtheorem{Lem}[Thm]{Lemma}
\newtheorem{Cor}[Thm]{Corollary}
\newtheorem{Prop}[Thm]{Proposition}
\newtheorem{Rm}{Remark}

\def\a{{\mathbb a}}
\def\C{{\mathbb C}}
\def\A{{\mathbb A}}
\def\B{{\mathbb B}}
\def\D{{\mathbb D}}
\def\E{{\mathbb E}}
\def\R{{\mathbb R}}
\def\P{{\mathbb P}}
\def\S{{\mathbb S}}
\def\Z{{\mathbb Z}}
\def\O{{\mathbb O}}
\def\H{{\mathbb H}}
\def\V{{\mathbb V}}
\def\Q{{\mathbb Q}}
\def\Cn{${\mathcal C}_n$}
\def\CM{\mathcal M}
\def\CG{\mathcal G}
\def\CH{\mathcal H}
\def\CT{\mathcal T}
\def\CF{\mathcal F}
\def\CA{\mathcal A}
\def\CB{\mathcal B}
\def\CD{\mathcal D}
\def\CP{\mathcal P}
\def\CS{\mathcal S}
\def\CZ{\mathcal Z}
\def\CE{\mathcal E}
\def\CL{\mathcal L}
\def\CV{\mathcal V}
\def\CW{\mathcal W}
\def\IC{\mathbb C}
\def\IF{\mathbb F}
\def\IK{\mathcal K}
\def\IL{\mathcal L}
\def\IP{\bf P}
\def\IR{\mathbb R}
\def\IZ{\mathbb Z}

\title{Isotoping  $2$-spheres in $4$-manifolds }
\author{Selman Akbulut}
\thanks{Partially supported by NSF grants DMS 0905917}
\keywords{}
\address{Department  of Mathematics, Michigan State University,  MI, 48824}
\email{akbulut@math.msu.edu }
\subjclass{58D27,  58A05, 57R65}
\date{\today}
\begin{abstract} 
Here we discuss an example of topologically isotopic but smoothly non-isotopic pair of $2$- spheres in a simply connected 4-manifold, which become smoothly isotopic after stabilizing by connected summing with $S^{2}\times S^{2}$.
\end{abstract}

\date{}
\maketitle

\setcounter{section}{-1}

\vspace{-.25in}

\section{The example}

in \cite{akmr} examples of topologically isotopic but smoothly nonisotopic spheres in simply connected $4$-manifolds, which become smoothly isotopic after connected summing with $S^{2}\times S^{2}$, were discussed. In this note we show that  such an example already follows from \cite{a2}.

\vspace{.1in}

\begin{figure}[ht]  \begin{center}
 \includegraphics[width=.41\textwidth]{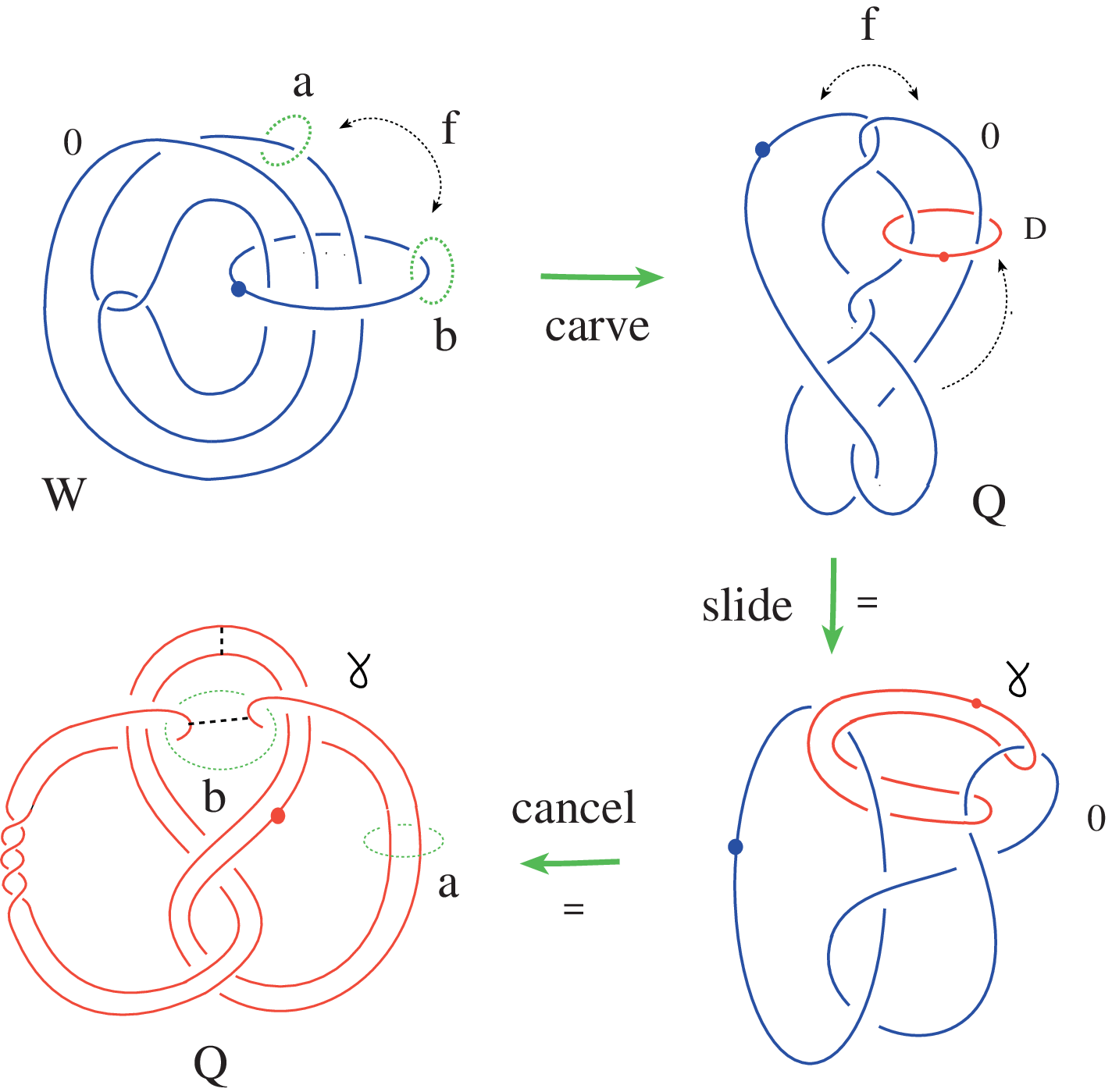}       
\caption{}      \label{s1} 
\end{center}
  \end{figure}

First review \cite{a2}: Let  $f:\partial W\to \partial W$ be the cork twisting involution of the Mazur cork $(W, f)$. Since  $f$ fixes the boundary $\partial D$ of a properly imbedded disk $D\subset W$  up to isotopy, $f(\partial D)$ bounds a disk in $W$ as well (isotopy in the collar union $D$), hence we can extend $f$ across the tubular neighborhood of $N(D)$ of $D$ by the carving process of  \cite{a1}. This results a manifold $Q = W- N(D)$ of Figure~\ref{s1},  homotopy equivalent to $B^{3}\times S^{1}$, and an involution on its boundary $\tau : \partial Q \to \partial Q$.

$\tau$ does not extend to $Q$ as a diffeomorphism (otherwise $f$ would extend to a self diffeomorphisim of $W$).  So  $\tau$ gives an exotic structure to $Q$ relative to its boundary (just as in  the cork case).  In \cite{a4} such $(Q,\tau)$'s are called anticorks because they live inside of corks $(W,f)$, and twisting $Q$ by the involution $\tau$ undoes the effect of twisting $W$ by $f$. Notice that the loop $\gamma= \partial D$ bounds two different disks in $B^4$ with the same complement $Q$ (where the identity map between their boundaries can not extend to a diffeomorphism inside), they are described by the two different ribbon moves indicated in the last picture of Figure~\ref{s1}. The two disks are the obvious disks which $\gamma$ bounds in the third picture of Figure~\ref{s1}, and the same disk after zero and dot exchanges of the figure.

\vspace{.1in}

\begin{figure}[ht]  \begin{center}
 \includegraphics[width=.53\textwidth]{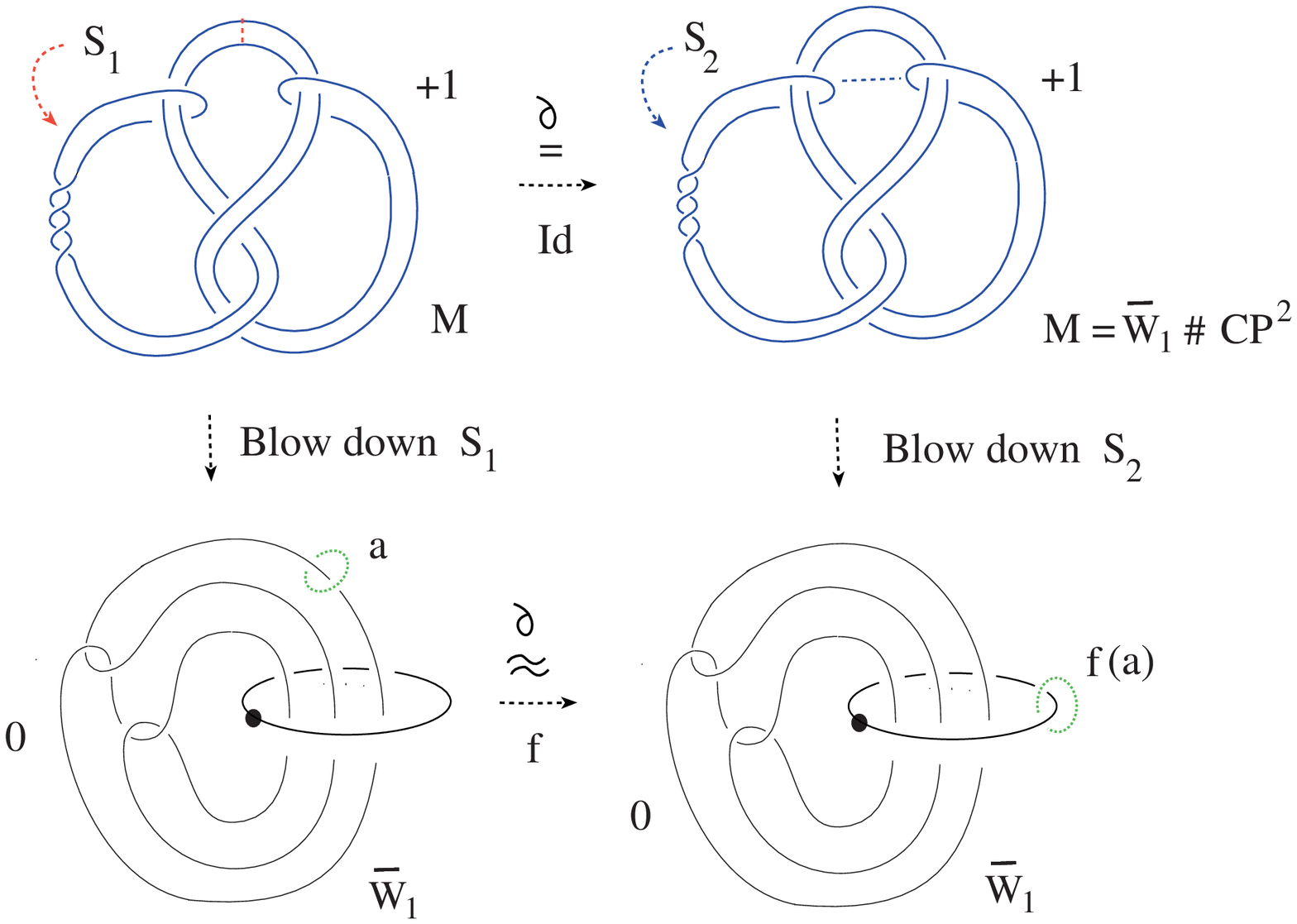}       
\caption{}      \label{s2} 
\end{center}
  \end{figure}

 Now let $M$ be the $4$-manifold obtained by attaching a $2$-handle to $B^4$ along the ribbon knot $\gamma $ of Figure~\ref{s1},  with $+1$ framing. Clearly $M$ has two imbedded $2$-spheres $S_{i}$, $i=1,2$ of self intersection $+1$ representing $H_{2}(M)\cong \Z$, corresponding to the two different $2$-disks which $\gamma$ bounds in $B^4$. Blowing down either $S_{1}$ or $S_{2}$ turns $M$ into the {\it positron} cork $\bar{W_{1}}$ of Figure~\ref{s2} (\cite{am}), and the two different blowings down process turns the identity map $\partial M \to \partial M $ to the cork involution $f:\partial \bar{W_{1}}\to \partial \bar{W_{1}}$, i.e. the maps in Figure~\ref{s2} commute (this can be seen by blowing down $\gamma$  of Figure~\ref{s1} by using the two different disks). Hence $S_{1}$ and $S_{2}$ are not smoothly isotopic in $M$ by any isotopy keeping $\partial M$ fixed, though they are topologically isotopic (by Freedman's theorem); but they are isotopic In $M\# S^{2}\times S^{2}$ rel boundary (since surgery corresponds turning the dotted circle to $0$-framed circle, in the third picture of Figure~\ref{s1}).

 Since  $M=\bar{W_{1}} \# \C\P^{2}$, this example shows that the blowing up $ \C\P^{2}$ operation undoes the cork twisting operation of  $(\bar{W_{1}}, f)$. Also since the Dolgachev surface $E(1)_{2,3}$ differs from its standard copy $\C\P^{2}\# 9\bar{{\C\P}^{2}}$ by twisting the positron cork  $(\bar{W_{1}}, f)$ inside (Theorem 1 of \cite{a3}), the manifold $M$ in this example can be made to be closed (without boundary).
 
 { \Rm Reader can check that the two ribbon disks of the ribbon knot in Figure~\ref{s1} are actually the same ribbon disks (isotoping the last picture of Figure~\ref{s1}, by forcing the two stands going through $b$ stay parallel, results the same picture except  the positions of $a$ and $b$ are exchanged) but $f$ induces nontrivial identifications on the boundaries of the ribbon complements.}

\end{document}